\documentclass[12pt,a4paper,oneside]{article}
\usepackage{amsmath,amstext,amssymb,amsopn,amsthm,mathrsfs}

\theoremstyle{plain}
\newtheorem{thr}{Theorem}[section]
\newtheorem{lem}[thr]{Lemma}
\newtheorem{prop}[thr]{Proposition}
\newtheorem{cor}[thr]{Corollary}
\newtheorem*{property*}{Property A}

\theoremstyle{definition}

\newtheorem{assump*}{Assumption}

\theoremstyle{remark}
\newtheorem{rem}[thr]{Remark}

\newcommand{\Rd}{\mathbb{R}^d}
\newcommand{\D}{\delta_D(x)\delta_D(y)}
\newcommand{\G}[2]{\widetilde{G}_D(#1,#2)}

\newcommand{\Z}{\int^{\infty}_{0}}
\newcommand{\ro}{\mathbf{r}_0}
\renewcommand{\leq}{\leqslant}
\renewcommand{\geq}{\geqslant}
\newcommand{\norm}[1]{|#1|}
\newcommand{\til}[1]{\widetilde{#1}}

\DeclareMathOperator{\diam}{diam}

\title{Estimates of Green Function for some perturbations of fractional Laplacian}
\author{ {\sc TOMASZ GRZYWNY and MICHA\L{} RYZNAR}\\
{\footnotesize Institute of Mathematics and Computer Science} \\
{\footnotesize Wroc\l{}aw University of Technology} \\
{\footnotesize Wyb{.} Wyspia\'nskiego 27, 50--370 Wroc\l{}aw, Poland}\\
{\footnotesize email: tomasz.grzywny@pwr.wroc.pl;
Michal.Ryznar@pwr.wroc.pl}}
\date{}

\begin{document}
\maketitle

\begin{abstract}
Suppose that  $Y(t)$ is a $d$-dimensional L\'{e}vy symmetric process for  which its L\'{e}vy measure differs from the L\'{e}vy measure of the isotropic $\alpha$-stable process ($0<\alpha<2$) by a finite signed measure. For a bounded Lipschitz set $D$ we  compare  the Green functions of the process $Y$ and its stable counterpart.
We prove a few comparability results either one sided or two sided. Assuming an additional condition about the difference of the densities of the L\'{e}vy measures, namely that it is of order of $|x|^{-d+\varrho}$ as $|x|\to 0$, where $\varrho>0$, we prove that the Green functions are comparable, provided $D$ is  connected.   

These results apply for example  to $\alpha$-stable relativistic process. This process was studied in  \cite {R, CS3}, where the bounds for its Green functions were proved for $d> \alpha$ and smooth sets.  In the paper we also considered one dimensional case for $\alpha\ge 1$ and proved that the Green functions for an open and bounded interval are comparable.
 \end{abstract}

\section{Introduction}

The  purpose of the paper is to study estimates of the Green
functions of bounded open sets of a symmetric L\'{e}vy process
$Y_t$, which lives on $\Rd$. We assume that its L\'{e}vy measure
is close in some sense, which we specify later,   to the L\'{e}vy
measure of the isotropic $\alpha$-stable process. From the point
of view of infinitesimal generators, the generator of the
semigroup corresponding to  $Y_t$ can be considered as a
perturbation of the fractional Laplacian by a bounded linear
operator.  The potential theory of the stable process was
extensively investigated in the recent years (see \cite{Bogdan1},
\cite{BogdanBycz1}, \cite{CS1}, \cite{K}) and the there are
several results providing the estimates of the Green functions of
$C^{1,1}$ bounded sets (see \cite{Kulczycki} and \cite{CS2}) or
even bounded Lipschitz sets (\cite{J}, \cite{Bogdan2}). We intend
to make a comparison of the Green function of the process  $Y_t$
and its stable counterpart. One of the first results in this
direction was contained in \cite{R}, where so called relativistic
$\alpha$-stable process was considered. This is a process which
characteristic function is of the form
\[E^{0}e^{i z\cdot Y_t}=e^{-t((|z|^2+m^{2/\alpha})^{\alpha/2}-m)},\quad z\in \Rd,\]
 where $0<\alpha<2$ and $m>0 $ is a
parameter. Observe that for $m=0$ it reduces to the isotropic
$\alpha$-stable process.  The main result of \cite{R} says that
the Green function of $C^{1,1}$ bounded set was comparable to the
Green function of the  isotropic $\alpha$-stable process if
$d>\alpha$. Later on that result was derived by a different method
in \cite{CS3}.  In the present paper we develop methods from
\cite{R} to derive several extensions of the results proved
therein. The main results are contained in the following two
theorems.

\begin{thr}\label{main1}Let $D\subset \Rd$ be a Lipschitz connected and bounded open set.
Suppose that $Y_t$ is a symmetric purely jump L\'{e}vy process in
$\Rd$ with $d\ge 1$ and $\nu^Y(x)$ is the density of its L\'{e}vy
measure. By $\til{\nu}(x)$ we denote the density of the L\'{e}vy
measure of the isotropic  stable process and by  $\til{G}_{D}$ its Green
function of $D$. Assume that
 $\sigma(x)=\til{\nu}(x)-\nu^Y(x)\ge 0, \ x\in\Rd$, and
 $\sigma(x)\leq c|x|^{\varrho-d}$ for
$|x|\leq 1$, where $c,\varrho>0$. Then there exists a
constant $C=C(d,\alpha,D,\varrho,c)$, such that
\[C^{-1}\til{G}_{D}(x,y)\leq G^Y_{D}(x,y)\leq C \til{G}_{D}(x,y),\]
for all $x,y\in D$.
\end{thr}

In the next theorem we remove  the assumption about positivity of the function
$\sigma$ at the cost of some mild assumption about the behaviour of  the density of the L\'{e}vy measure.

\begin{thr}\label{main2}With the same notation as in the previous theorem assume that
there are positive constants $c$ and $\varrho$ such that
 $|\sigma(x)|\leq c |x|^{-d+\varrho}$ for $|x|\leq
1$,  and  $\nu^Y(x)$ is bounded on $B^c(0,1)$. Then there is a
constant $C=C(d,\alpha, D, \varrho, \sigma)$ such that for any
$x,y\in D$,
\[C^{-1}\til{G}_{D}(x,y)\leq G^Y_{D}(x,y)\leq C \til{G}_{D}(x,y).\]
\end{thr}

Observe that in the first theorem the assumption about the
positivity of $\sigma$ enables us  not to assume anything about
the behaviour of $\nu^Y(x)$ away from the origin except it has to
be dominated by  $\til{\nu}$. For example $\nu^Y(x)$ can vanish
outside some neighborhood  of the origin. Of course that
assumptions  are readily checked  for the relativistic process
(see \cite{R} for the description of the L\'{e}vy measure), so the
theorem extends to Lipschitz bounded domains the main result of
\cite{R} (see also \cite{CS3}). In addition,  note that it covers
the one-dimensional case for $\alpha\ge 1$, which was not treated
in the neither papers cited above. Actually both papers assumed
$d\ge 2$ but the proofs remain valid for $d> \alpha$.  To our best
knowledge  the one dimensional result is a new one which fills the
gap in the potential theory of the relativistic process.

The methods we apply are elementary and are based on the fact that
for any two pure jump processes  such that the difference of their
L\'{e}vy measures is a {\em positive} and  {\em finite} measure
one can represent one of the processes as a sum of the other and
an independent compound Poisson  process. A  different approach in
taken in \cite{CS3}, where the problem in $C^{1,1}$ case  was
tackled by so called drift transform technique. After obtaining
the main results of the present paper the authors found on the
website of Panki Kim a paper of Kim and Lee \cite{KimLee} with similar results
as ours but even for more general sets (so called $\kappa$-fat
sets). The method they use is essentially designed in \cite{CS3},
so our methods and  results can be viewed as an alternative
approach to the problem of comparing the Green functions. Moreover
our method can handle the situation when a L\'{e}vy measure
vanishes outside  some neighborhood  of the origin which seems not
be an option in the other method used in \cite{CS3} or
\cite{KimLee}.

The paper is organized in the following way. In Section 2 we set
up the notation and provide necessary definitions and basic facts
needed in the sequel. At first we do not assume that $Y_t$ is
compared with the stable process but we sometimes  work in
slightly more general setup. Namely some of the results are
formulated in such a  way that  $Y_t$ is compared with another
L\'{e}vy process $X_t $ under the appropriate assumptions about
their L\'{e}vy measures.  In Section 3 we prove the main estimates
along with some other related results. To prove Theorem
\ref{main2}  we first prove the estimates for sets of small
diameter and then use it to prove Boundary Harnack Principle (BHP) for
the process $Y_t$ in the case when its  L\'{e}vy measure dominates
the L\'{e}vy measure of the isotropic $\alpha$-stable process.

\section{Preliminaries}
In  $\Rd$, $d\geq 1$, we consider a symmetric L\'{e}vy processes
$X_t$ such that its characteristic triplet is equal to
$(0,\nu,0)$, where $\nu$ is its (nonzero) L\'{e}vy measure. That
is its characteristic function is given by
\[E^{0}e^{i z\cdot X_t}=e^{-t\int_{\Rd}(1-\cos(z\cdot w))\nu(dw)},\quad z\in {\bf R}^d.\]
If the measure $\nu$ is  absolutely continuous with respect to
the Lebesque measure then by $\nu(x)$ we denote its density. By
$p(t,x,y)$ we denote the transition densities of $X_t$, which are
assumed to be bounded and  defined for every $x,y\in \Rd $. The
potential kernel for $X_t$ is given by
\[U(x,y)=U(x-y)=\Z p(t,x-y)dt.\]

We use the notation $C=C(\alpha,\beta,\gamma,\dots)$ to denote that the
constant $C$ depends on $\alpha,\beta,\gamma,\dots$. Usually  values of
constants may change from line to line, but they are always  strictly
positive and finite. Sometimes we skip in notation that constants
depend on usual quantities (e.g. $d,\alpha$). Next, we give some
definitions. We use $f\approx g$ on $D$ to denote that the functions $f$ and $g$
are comparable, that is there exists a constant $C$ such that
$$C^{-1}f(x)\le g(x) \le C f(x), \quad x\in D.$$
Let $D\subset \Rd $ be an  open set.
By $\tau_D$ we denote the first exit time from $D$ that is
$$\tau_D=\inf\{t>0: X_t\notin D\}.$$

Next, we  investigate  boundness of the first moment of $\tau_D$.
\begin{lem}\label{momentzbogr:this}
For any bounded open set $D$ there exists a constant $C=C(D)$ such
that
\[\sup_{x\in \mathbb{R}^d}E^x\tau_D \leq C.\]
\end{lem}
\begin{proof}
The proof of this lemma follows by the same arguments as in the
classical case for the Brownian motion (see \cite{CZ}). The
argument therein requires the existence of $t_0>0$ such that
$\sup_{x\in\Rd}P^x(X_{t_0}\in D)< 1$. However, repeating the steps
from Lemma 48.3 in \cite{Sato},   one can obtain that
\[\sup_{x\in\Rd}P^x(X_t\in D)=O(t^{-1/2}),\quad  t\to \infty .\]
\end{proof}

In order to study the killed process on exiting of $D$ we construct
its transition densities by the classical formula
$$p_D(t,x,y)=p(t,x,y)-r_D(t,x,y),$$
where
$$r_D(t,x,y)=E^x[t\geq \tau_D;p(t-\tau_D,X_{\tau_D},y)].$$
The arguments used for Brownian motion (see eg.  \cite{CZ}) will
prevail in our case and one can easily show that $p_D(t,x,y), \
t\ge 0,$ satisfy the Chapman-Kolmogorov equation (semigroup
property). Moreover the transition density $p_D(t,x,y)$ is a
symmetric function $(x,y)$ a.s. Assuming some other mild
conditions on the transition densities of the (free) process one
can actually show that  $p_D(t,x,y)$  can be chosen as continuous
functions of $(x,y)$. Next, we define the Green function of the
set $D$,
$$G_D(x,y)=\Z p_D(t,x,y)dt.$$
Let us  see that the integral is well defined, because
$$\int_D G_D(x,y)dy=\int_D\Z p_D(t,x,y)dtdy=\Z P^x(\tau_D>t)dt=E^x\tau_D<\infty.$$
Hence for every $x\in \Rd$ the Green function $G_D(x,y)$ is well defined $(y)$ a.s.
Again under the assumptions which make
$p_D(t,x,y), \ t>0,$ continuous functions in arguments $x,y$ one can show that
the Green function is a continuous (in extended sense) function on $D\times D$.

It is well known that if the L\'{e}vy measure is  absolutely
continuous  with respect to the Lebesque measure then  the
distribution of $X_{\tau_D}$ restricted to $\overline{D}^c$ is
absolutely continuous as well (see Ikeda Watanabe) and the density
is given by so called Ikeda-Watanabe formula:
$$P_D(x,z)= \int_D G_D(x,y)\nu(y-z)dy, \quad (x,z)\in D\times \overline{D}^c $$
We call $P_D(x,z)$ the Poisson kernel. Under some other mild
conditions $X_{\tau_D}$  has zero probability of belonging to the
boundary od $D$ so in this case  the Poisson kernel fully describe
the distribution of the exiting point.

We say that measurable function $u$ is {\em harmonic} with respect
to $X_t$ in an open set $D$ if for every bounded open set $U$
satisfying $\overline{U}\subset D$,
$$ u(x) = E^x u(X_{\tau_U}),\quad x\in U.$$
Whereas if
$$ u(x) = E^x u(X_{\tau_D}),\quad x\in D,$$
then we say that $u$ is {\em regular harmonic} with respect to
$X_t$ in an open set $D$.

The following lemma is a simple  consequence of Lemma
\ref{momentzbogr:this} and  boundness of $p(t,x)$.
\begin{lem}\label{gestzabity2:this}
For any $x\in D$ and  $t\ge 1$ we have
$$p_D(t,x,y)\leq C(X) \frac{E^x\tau_D E^y\tau_D}{t^2}\quad (y) \text{ a.s. } .$$
\end{lem}

\begin{proof}
 Observe that for $s\geq 0,$
\[\sup_{x,y\in D}p_D(s+1/2,x,y)\leq \sup_{x,y\in\Rd}p(s+1/2,x-y)=\sup_{x\in\Rd}p(1/2,\cdot)*p(s,x)
\leq \sup_{x\in\Rd}p(1/2,x)=C_1.\]
Hence, by the Chapman-Kolmogorov equation we obtain for $t\geq 1$
and $(y)$ a.s.
$$p_D(t,x,y)=\int_D p_D(t/2,x,z)p_D(t/2,z,y)dz \leq C_1 P^x(\tau_D>t/2).$$
Applying again the Chapman-Kolmogorov equation together with the above inequality we get
$$p_D(t,x,y)\leq C_1 P^x(\tau_D>t/4)\int_D p_D(t/2,z,y)dz=C_1P^x(\tau_D>t/4)P^y(\widehat{\tau}_D>t/2),$$
where $\widehat{\tau}_D=\inf\{t>0:-X_t\in D\}$. But the process $X_t$
is symmetric, so $\{X_t\}\stackrel{D}{=}\{-X_t\}$. Hence
$$P^y(\widehat{\tau}_D>t/2) = P^y(\tau_D>t/2).$$
Therefore, we have
$$p_D(t,x,y)\leq C_1 P^x(\tau_D>t/4) P^y(\tau_D>t/2).$$
The application of Chebyshev's inequality completes the proof.
\end{proof}
\begin{rem}\label{gestzabityStable}
If $X_t$ is isotropic stable process then by similar arguments we have for $t>0$ and
$x,y\in D$,
$$p_D(t,x,y)\leq C(\alpha,d) \frac{E^x\tau_D E^y\tau_D}{t^{2+d/\alpha}}.$$
\end{rem}

In one of our general results (Theorem \ref{general}) we require
the following property which exhibits a relation between moments
of the exiting times and the Green function.

\begin{property*}Suppose that there is a constant $ c=c(D)$ such
that
\begin{equation*}
E^x\tau_D E^y\tau_D\le c G_D(x,y), \quad x,y\in D.
\end{equation*}
\end{property*}

At the first glance the above condition looks a bit restrictive but
actually it holds  in the   stable case (\cite{K}, \cite{CS1},
\cite{Banuelos}) and usually it is derived as a consequence of the
intrinsic ultracontractivity of the killed process. In the recent
paper of the first author (see \cite{GR}) the
intrinsic ultracontractivity is studied under much
broader assumptions. For example the above property holds if
$p_D(t,\cdot,\cdot)$ is continuous in $x,y$ and the Lebesgue measure is
absolute continuous with respect to the L\'{e}vy measure.

From now we  consider two symmetric L\'{e}vy processes $Y_t$ and
$X_t$ such that a signed measure $\sigma=\nu^X-\nu^Y$ is finite,
where $\nu^Y$, $\nu^X$ are L\'{e}vy measures of $Y_t$ and $X_t$
respectively. We  use that notational convention throughout the
whole paper, e.g. we denote the transition density of $X_t$ by
$p^X(t,x)$ and the transition density of $Y_t$ by $p^Y(t,x)$.
Later on we specify one of the processes, say  $X_t$,  to be the
isotropic stable process. The aim of this paper is to provide some
comparisons between the two process in various aspects of which
the relationship of the Green functions is our main target. Some
of the results are general but our typical situation is a
comparison between the isotropic stable process and another
process with the L\'{e}vy measures sufficiently close to each
other.

With the  assumption that $\sigma=\nu^X-\nu^Y$ is finite we can
write the  following formula comparing infinitesimal generators on
$L^1(\Rd)$ of these processes

\begin{equation*}
\mathcal{A}^Y=\mathcal{A}^X-P,\quad \text{ where } P\varphi(x)= \sigma
* \varphi (x) - \sigma(\Rd) \varphi(x).
\end{equation*}
The fact that $P$ is a bounded operator implies that the domains
of these generators coincide.

As mentioned above, very often the process ${X}_t$ is taken to be
the isotropic  $\alpha$-stable process, $0<\alpha<2$ . To emphasize
its role we denote it  by $\til{X}_t$.  That  process has the
following characteristic function:
\[E^{0}e^{i z\cdot \til{X}_t}=e^{-t|z|^{\alpha}},\quad z\in \Rd.\]
 From now on, we will use the  tilde sign to denote
functions, measures and etc. corresponding to $\til{X}_t$. For
example its L\'{e}vy measure is given by the formula
$$\til{\nu}(B)=\int_B
\mathscr{A}(-\alpha,d)|x|^{-d-\alpha}dx,$$ where
$\mathscr{A}(\rho,d)=\frac{\Gamma((d-\rho)/2)}{\pi^{d/2}2^{\rho}|\Gamma(\rho/2)|}$.
The potential kernel which is well defined for $\alpha<d$ is given
by
$$\til{U}(x)= \mathscr{A}(\alpha,d)|x|^{\alpha-d},\quad x\in \Rd.$$


The next two lemmas provide basic tools for examining the
relationship between the Green functions. In the first we compare
the moments of exiting times only under the assumption that
$\sigma=\nu^X-\nu^Y$ is a finite signed measure, while in the second
we require that $\sigma$ is nonnegative.  This assumption
provides us with a nice inequality involving the transitions
densities. However the both lemmas  already appeared in \cite{R}
under some additional assumptions, we deliver the proofs for the
reader convenience.
\begin{lem}\label{momentyogr:this}
Let $D$ be a bounded open set and $\sigma=\nu^X-\nu^Y$ be finite.
Then we have on $D,$
\[ E^x \tau^X_D\approx  E^x\tau^Y_D.\]
\end{lem}

\begin{proof} Suppose that the Jordan decomposition of $\sigma=\sigma_{+}-\sigma_{-}$.
Let $V_t$ be a compound Poisson process independent of $X_t$ with the
L\'{e}vy measure $\sigma_{-}$ and $V^{'}_t$ be a compound Poisson
process  independent of $Y_t$ with the L\'{e}vy measure $\sigma_{+}$.
We put $Z_t=X_t+V_t$, then of course we have
$\{Z_t\}\stackrel{D}{=}\{Y_t+V^{'}_t\}$. Hence, it's enough to
show that $E^x\tau^Z_D\approx E^x\tau^X_D$.

Let us define a stopping time $T$ by $T=\inf\{t>0: V_t\neq 0 \}$.
The processes $X_t$ and $V_t$ are mutually independent, therefore
$X_t$ and $T$ are independent as well. Besides, $Z_t=X_t$ for
$0\leq t < T$. We set  $m=\sigma_{-}(\Rd)$.

First, we claim that $ E^x (\tau^X_D) \leq 2 E^x (\tau^X_D\wedge t)$
for $t$ large enough. Indeed, by the Markov Property and
Lemma \ref{momentzbogr:this} we have
\begin{eqnarray*}
E^x\tau^X_D &=& E^x (\tau^X_D\wedge t) + E^x({\tau^X_D>t};\tau^X_D-t) =
E^x (\tau^X_D\wedge t)+ E^x({\tau^X_D>t};E^{X_t}\tau^X_D)\\
&\leq& E^x (\tau^X_D\wedge t) + C P^x(\tau^X_D>t) \leq E^x
(\tau^X_D\wedge t) + C \frac{E^x\tau^X_D}{t},
\end{eqnarray*}
which proves our claim for $t\geq 2C$.

Because $\tau^Z_D\wedge T=\tau^X_D\wedge T$, so by independence $T$
and $X_t$ we get
\begin{eqnarray*}E^x \tau^Z_D &\geq& E^x (\tau^Z_D\wedge T) =
E^x (\tau^X_D\wedge T) =
 \int^{\infty}_0 E^x(\tau^X_D\wedge t)me^{-mt}dt \\
 &\geq& \int^{\infty}_{2C} E^x(\tau^X_D\wedge t)me^{-mt}dt\geq \frac{1}{2}e^{-2 C m }
 E^x\tau^X_D.
 \end{eqnarray*}

Now, we prove the upper bound
\begin{eqnarray*}
E^x \tau^Z_D &=&  E^x (\tau^Z_D\wedge T) +
E^x({\tau^Z_D>T};\tau^Z_D-T)\\
&\leq & E^x \tau^X_D + E^x({\tau^Z_D>T};E^{Z_T}\tau^Z_D)\\
&\leq & E^x \tau^X_D + C P^x(\tau^Z_D>T),
\end{eqnarray*}
but \[P^x(\tau^Z_D>T)\leq P^x(\tau^X_D\geq T)=m \Z P^x(\tau^X_D\geq
t)e^{-mt}dt\leq m E^x \tau^X_D,\]
which ends the proof.
\end{proof}

\begin{lem}\label{gestzabity1:this}
Suppose that $\sigma=\nu^X-\nu^Y$ is a nonnegative finite measure
and $D$ is an open set. Then for any $x\in D$ and $t>0$,
\[
p^Y_D(t,x,\cdot) \leq  e^{mt} p^X_D(t,x,\cdot)\quad a.s.\text{
}.\] If, in addition, we assume that $p^Y(t,\cdot)$ and
$p^X(t,\cdot)$  are continuous then we have for $x,y\in D$,
\[r^Y_D(t,x,y) \leq  e^{2mt} r^X_D(t,x,y).\]

\end{lem}
\begin{proof}
We put $m=\sigma(\Rd)<\infty$,
 and  define a compound Poisson process $V_t$
with the L\'{e}vy measure $\sigma$ independent of $Y_t$. A random
variable \begin{equation}\label{defT}T=\inf\{t\geq 0: V_t\neq
0\}\end{equation}  has the exponential distribution with intensity
$m$. Then $Y_t$ and $T$ are independent and for $0\leq t< T$ we
have $X_t=Y_t$.

Let $A$ be a Borel subset of $D$. Since $Y_t=X_t$, for $t<T$ we infer that
$\{\tau^Y_D>t\}\cap\{T>t\}=\{\tau^X_D>t\}\cap\{T>t\}$. \mbox{By
independence of} $Y_t$ and $T$
\begin{eqnarray*}
P^x(t<\tau^Y_D;Y_t\in A)P^x(T>t)&=&P^x(t<\tau^Y_D;Y_t\in A;T>t)\\
&=&P^x(t<\tau^X_D;X_t\in A;T>t)\\
&\leq& P^x(t<\tau^X_D;X_t\in A).
\end{eqnarray*}
So we obtain  that $(y)$ a.s. ,
\[p^Y_D(t,x,y)P^x(T>t)\leq p^X_D(t,x,y).\]
But $T$ has the exponential distribution with intensity $m$, that
is $P^x(T>t)=e^{-mt}$.

The second  inequality is proved analogously,  using the first
with $D=\Rd$ in the intermediate step. Moreover the continuity of
$p^Y(t,\cdot)$ and $p^X(t,\cdot)$ is required to justify the last
step:
\begin{eqnarray*}
r^Y_D(t,x,y)e^{-mt}&=&E^x[t\geq\tau^Y_D;p^Y(t-\tau^Y_D,Y_{\tau^Y_D},y)]P^x(T>t)\\
&=&E^x[\tau^Y_D\leq t<T;p^Y(t-\tau^Y_D,Y_{\tau^Y_D},y)]\\
&=&E^x[\tau^X_D\leq t<T;p^Y(t-\tau^X_D,X_{\tau^X_D},y)]\\
&\leq&e^{mt}E^x[\tau^X_D\leq t;p^X(t-\tau^X_D,X_{\tau^X_D},y)]\\
&=&e^{mt}r^X_D(t,x,y).
\end{eqnarray*}
\end{proof}

The next lemma is a sort of a comparison between transition
densities in that sense that a "nice" behaviour of them for one
process implies that the  transition densities of the second are
uniformly bounded away from zero. The "nice" behaviour for example
is present if the first process is the isotropic stable process.
We use that result in the sequel to assure that the transition
densities of the killed process are continuous  and to assure the
property A. We define an exponent of a signed finite measure
$\sigma$ by
$$\exp\{\sigma\}(A)=e^{-\sigma(\Rd)}\sum^{\infty}_{n=0}\frac{\sigma^{*n}(A)}{n!},
\quad \text{where } A\subset\Rd \text{ is a Borel set}.$$

\begin{lem}\label{BoundDensity}
Suppose that $\nu^X$ and $\nu^Y$ are absolutely continuous and
$\,\sigma(x)=\nu^X(x)-\nu^Y(x)$ is an integrable function
such that $|p^X(t,\cdot)*\sigma(x)|+|\sigma(x)|\le c_1$ for $|x|\ge \delta$
and $t\le 1$. If $p^X(t,x)\leq c_2t^{-\zeta}$ for $t\le 1$, where $\zeta>0$,
and $p^X(t,x)\leq c_3(\delta)$ for $|x|\ge \delta$,
then there is a constant $C$ such that
$$p^ Y(t,x)\le C,  \quad |x|\ge \left(\left[\zeta\right]\vee1\right)
\delta \text{ and } t>0. $$
\end{lem}

\begin{proof}
 Suppose that $\int_{\Rd}|\sigma(x)|dx = M < \infty$.
 We put $\int_{\Rd}\sigma(x)dx=m$. We can write
$$  {p}^Y(t,x) = p^X(t,\cdot)*\exp\{-t \sigma\}=
 p^X(t,x)e^{tm}+\sum_{n=1}^\infty  \frac  {(-t)^n p^X(t,\cdot)*\sigma^{*n} (x)}{n!}e^{tm}.$$

 Observe that $|p^X(t,\cdot)*\sigma^{*n} (x)|\le \sup_{y\in R^d}p^X(t,y)M^n\leq
 c_2 \frac {M^n}{t^{\zeta}}$, so for $t\le 1$ we have
\begin{equation}\label{ogon_gest}
|\sum_{n\ge \zeta}^\infty  \frac{(-t)^n p^X(t,\cdot)*\sigma^{*n}
(x)}{n!}e^{tm}|\le C\sum_{n\ge \zeta}^\infty
\frac{t^{n-\zeta}M^n}{n!}=Ce^{M}<\infty.
\end{equation}
Now, we show that if $|p^X(t,\cdot)*\sigma(x)|+|\sigma(x)| \le c(1)$ for
$|x|\ge \delta$ and $t\le 1 $ then
\begin{equation}\label{pocz_gest} |p^X(t,\cdot)*\sigma^{*n}(x)|\le
c(n), \quad |x|\ge n \delta.\end{equation}  We assume (\ref{pocz_gest}) for
$n$ and we prove it for $n+1$. Observe that
\begin{eqnarray*}
|p^X(t,\cdot)*\sigma^{*n+1}(x)|&\leq&
\int_{B^c(x,n\delta)}|p^X(t,\cdot)*\sigma^{*n}(x-y)||\sigma(y)|dy
+\\
&&+\int_{B(x,n\delta)}|p^X(t,\cdot)*\sigma^{*n}(x-y)||\sigma(y)|dy\\
&\le& c(n)M + c_1 M^n,
\end{eqnarray*}
because if $y\in B(x,n\delta)$ then $|y|\ge |x|-|x-y|\ge \delta$.
Combining (\ref{ogon_gest}) and (\ref{pocz_gest}) and using that
$p^X(t,x)\le c(\delta)$ for $|x|\ge \delta$ we end the proof for $t\leq 1$.

Next,  for $t>1$ we have
\[\sup_{x\in\Rd}p^Y(t,x)= \sup_{x\in\Rd}p^Y(1,\cdot)*p^Y(t-1,x)\leq  \sup_{x\in\Rd}p^Y(1,x)=C,\]
which proves the conclusion for $t>1$.
\end{proof}
The following lemma is an attempt to find a condition under which
the potential kernel of a process is comparable at the vicinity of
the origin  with the stable potential kernel. It will play an
important role in proving the upper bound for the Green function
$G_D^Y$ by its stable counterpart (see Theorem \ref{GreenBig}).

\begin{lem}\label{potentialcomp}
Let $d>\alpha$. Let $-\sigma= \nu^Y -\tilde \nu$ be a nonnegative
finite measure such that  $ \tilde{U}*(-\sigma) (x)\le C
\tilde{U}(x)$ for $|x|\le 1$ then for some constant $C>1$,
  $$C^{-1}\tilde{U}(x)\le U^Y(x)\le C\tilde{U}(x), \quad |x|\le 1.$$
\end{lem}

\begin{proof}

 Suppose that $-\sigma = \nu^Y - \tilde{\nu}\ge 0$. Let $-\sigma(\Rd)=m>0$. We can write
$$  {p}^Y(t,x) = \widetilde{p}(t,\cdot)*\exp\{-t \sigma\}=
 \widetilde{p}(t,x)e^{-tm}+\sum_{n=1}^\infty  \frac{t^n \widetilde{p}(t,\cdot)*(-\sigma)^{*n} (x)}{n!}e^{-tm}.$$

 Observe that $\widetilde{p}(t,\cdot)*(-\sigma)^{*n} (x)\le
 \sup_{y\in \Rd}\widetilde{p}(t,y)m^n= C \frac {m^n}{t^{d/\alpha}}$ so for $n> d/\alpha -1$ we have
\begin{eqnarray*}
\int_0^\infty \frac  {t^n \widetilde{p}(t,\cdot)*(-\sigma)^{*n}
(x)}{n!}e^{-tm}dt&\le&  C
\int_0^\infty \frac  {t^{n-d/\alpha}m^n }{ n!}e^{-tm}dt\\
 &\le& C\frac  {\Gamma(n+1-d/\alpha) }{ n!} m^{d/\alpha +1}\le C \frac{m^{d/\alpha +1}}{n^{d/\alpha}}.
\end{eqnarray*}
This implies that
\begin{equation}\label{ogon}
\int_0^\infty\sum_{n>d/\alpha -1}^\infty  \frac  {t^n
\widetilde{p}(t,\cdot)*(-\sigma)^{*n} (x)}{n!}e^{-tm}\le
C\sum_{n>d/\alpha -1}^\infty \frac   {m^{d/\alpha
+1}}{n^{d/\alpha}}=c(\alpha,m,d)<\infty.
\end{equation}

Next estimating $t^n e^{-tm}\le C(n,m)<\infty$ we have

$$\int_0^\infty \frac  {t^n \widetilde{p}(t,\cdot)*(-\sigma)^{*n} (x)}{n!}e^{-tm}dt\le
C(n,m)\,\til{U} *(-\sigma)^{*n} (x).$$ Let $\tilde{U}(x)=
\frac{\mathscr{A}}{|x|^{d-\alpha}}$. If we assume that $
\tilde{U}*(-\sigma) (x)\le C \tilde{U}(x)$ for $|x|\le 1$ then we claim that
\begin{equation}\label{iteracja}
\tilde{U}*(-\sigma)^{*n} (x)\le C(n)\tilde{U}(x), \quad |x|\le 1.
\end{equation}
We check this for $n=2$ since the general case will follow by induction.

\begin{eqnarray*}
\tilde{U}*\sigma^{*2} (x)&=& \int_{B(x,1)} \tilde{U}*(-\sigma)
(x-y)(-\sigma)(dy) + \int_{B^c(x,1)}
\tilde{U}(x-y)\sigma^{*2}(dy)\\
&\le& C\int_{B(x,1)} \tilde{U}(x-y)(-\sigma)(dy)+ \mathscr{A}\,m^2\\
&\le& C^2\tilde{U}(x) + \mathscr{A}\,m^2 \le  C(2) \tilde{U}(x),
\end{eqnarray*}
because  $lim_{|x|\to 0}\tilde{U}(x)=\infty$. By (\ref{ogon}) and
(\ref{iteracja}) we conclude that $ U^Y(x)\le C \tilde{U}(x), \
|x|\le 1$.

Getting the reverse inequality is almost immediate since
$\widetilde{p}(t,x)\le e^{tm} {p}^Y(t,x)$ (Lemma
\ref{gestzabity1:this} with the fact that $\til{p}(t,\cdot)$ and
$p^Y(t,\cdot)$ are continuous).  The following estimate is well known:
\begin{equation}\label{ogrGestStab}
\til{p}(t,x)\leq C(d,\alpha)\,\left( t^{-d/\alpha}\wedge
\frac{t}{|x|^{d+\alpha}}\right).
\end{equation}
Hence for $|x|\leq 1$,
$$\tilde{U}(x)\leq C \int_0^1\widetilde{p}(t,x)dt,$$ for some constant $C=C(d,\alpha)$.
Therefore
$$ \tilde{U}(x)\leq \int_0^1\widetilde{p}(t,x)dt \le
 e^{m}\int_0^1 {p}^Y(t,x)dt \le e^{m}{U}^Y(x),$$
for $\ |x|\le 1$.
\end{proof}

\begin{rem}If $-\sigma(x)$ is a nonnegative density of a finite measure and

$$-\sigma(x)\le C |x|^{-d+\varrho}, \quad |x|\le 1,$$

where $\varrho > 0$ then the condition $ \tilde{U}*(-\sigma)
(x)\le C \tilde{U}(x)$  for $|x|\le 1$ is satisfied.
\end{rem}

The last lemma in this section is intended to treat the
one-dimensional recurrent case while comparing two processes of
which one is a stable one. This case is different from the
transient one and requires somewhat different arguments.

\begin{lem}\label{potlematuzup:this}
Let $d=1$, $\alpha\ge 1$ and $0<t_0\le 1$. Suppose that $\sigma =
\til{\nu}-\nu^Y$ is a finite measure. Then there exists a constant
$C=C(m,M)$ such that
\[\int^{t_0}_0 |\widetilde{p}(t,x)-e^{-2mt}p^Y(t,x)|dt\leq C\,t^{2-1/\alpha}_0, \]
where $m=\sigma(\mathbb{R})$ and $M=|\sigma|(\mathbb{R})$.
\end{lem}

\begin{proof}
Let $\sigma(\mathbb{R})=m$ and $|\sigma|(\mathbb{R})=M>0$. We can
write
$$  {p}^Y(t,x) = \widetilde{p}(t,\cdot)*\exp\{-t \sigma\}=
 \widetilde{p}(t,x)e^{tm}+\sum_{n=1}^\infty  \frac  {(-t)^n \widetilde{p}(t,\cdot)*\sigma^{*n} (x)}{n!}e^{tm}.$$
Next $|\widetilde{p}(t,\cdot)*\sigma^{*n} (x)|\le \sup_{y\in
\mathbb{R}}\widetilde{p}(t,y)M^n= C \frac {M^n}{t^{1/\alpha}}$.
Using this estimate we obtain

\begin{eqnarray*} |\widetilde{p}(t,x)-e^{-2mt}p^Y(t,x)| &=& \left|\widetilde{p}(t,x)(1-e^{-mt})
-\sum_{n=1}^\infty \frac  {(-t)^n \widetilde{p}(t,\cdot)*\sigma^{*n} (x)}{n!}e^{-mt}\right|\\
&\le& \widetilde{p}(t,x)(1-e^{-mt})+ \frac
{C}{t^{1/\alpha}}\sum_{n=1}^\infty \frac  {(tM)^n} {n!}e^{-mt}
\end{eqnarray*}
From the above   it easily follows that there is a constant $C
=C(m, M)$ such that

 $$|\widetilde{p}(t,x)-e^{-2mt}p^Y(t,x)|\le  C t^{1-1/\alpha}, \quad t\le 1.$$
Now the conclusion follows by integration.
\end{proof}

\section{Comparability of the Green functions}
In this section we prove our main results. We start with a general
one-sided estimate of Green functions.
 \begin{thr}\label{general} Let $D$ be a bounded open set and
 a finite measure $\sigma=\nu^X-\nu^Y$ be nonnegative.
 Suppose that for one of the processes $X_t$ or $Y_t$
 its  Green function satisfies the property A.
  Then there exists a constant
$C=C(\sigma,D,\alpha,d)$ such that for $x\in D$,
\[G^Y_D(x,y)\leq C G^X_D(x,y)\quad (y)\text{ }a.s. .\]
\end{thr}
\begin{proof}
 Denote $\sigma(\Rd)=m$.
From Lemmas \ref{gestzabity2:this} and \ref{gestzabity1:this} we
get $(y)$ almost surely
\begin{eqnarray*}G^Y_D(x,y) &=& \int^{t_0}_0
p^Y_D(t,x,y)dt+\int^{\infty}_{t_0} p^Y_D(t,x,y)dt\\&\leq& e^{mt_0}
\int^{t_0}_0 p^X_D(t,x,y)dt+C_1\int^{\infty}_{t_0}t^{-2} E^{x}\tau^Y_D
E^{y}\tau^Y_D dt,
\end{eqnarray*}
for $t_0\ge 1$. Hence
\[G^Y_D(x,y) \leq cG^X_D(x,y)+\frac{C_1}{t_0} E^{x}\tau^Y_D E^{y}\tau^Y_D.\]
If  $Y_t$ satisfies
\begin{equation}\label{Greenogrdol}E^{x}\tau^Y_DE^{y}\tau^Y_D\le C_2 G^Y_D(x,y),\end{equation}
then for $t_0=\max\{1,2C_1C_2\}$ we get
\[G^Y_D(x,y) \leq 2cG^X_D(x,y).\]
Now,  suppose that (\ref{Greenogrdol}) holds for $X_t$. Then
by Lemma \ref{momentyogr:this} we have
\[G^Y_D(x,y) \leq cG^X_D(x,y)+C_3 E^{x}\tau^X_D E^{y}\tau^X_D\leq C G^X_D(x,y),\]
which ends the proof.
\end{proof}

Kulczycki in \cite{K} showed that for the isotropic $\alpha$-stable
process the property A  is satisfied for any bounded open set $D$, so
we obtain the following.
\begin{cor}\label{stableGreenUp}
Let $D$ be a bounded open set. If $\sigma=\til{\nu}-\nu^Y$ is  a nonnegative and finite measure
 then there is a constant $C$ such that
\[G^Y_D(x,y)\leq C \til{G}_D(x,y).\]
 If $\nu^Y-\til{\nu}$ is a nonnegative and finite measure  then
\[\til{G}_D(x,y)\leq C G^Y_D(x,y).\]
\end{cor}

Suppose that $p^X_D(t,x,\cdot)$ and $p^X_D(t,\cdot,x)$ are
continuous for any $x\in D$. If the Lebesgue measure is absolutely
continuous with respect to the L\'{e}vy measure of $X_t$, then the
following theorem is true for any bounded open set $D$. Whereas if
there exists a radius $r>0$ such that density $\nu^X_{ac}$ of the
absolute continuous part of the L\'{e}vy measure satisfies
$$\inf_{x\in B(0,r)}\nu^X_{ac}(x)>0,$$
then the following theorem holds for any bounded and connected
Lipschitz domain $D$ (see \cite{GR}).
\begin{thr}\label{IUThr:this}
For every $t>0$ there is a constant $c=c(t,D,\alpha)$ such that
\[c E^x \tau^X_D E^y\tau^X_D\leq p^X_D(t,x,y),  \qquad x,y\in D.\]
\end{thr}
If we integrate the above inequality with respect to $dt$ we get the  property A for $X_t$
\begin{equation*}\label{intrinsic}
C E^x\tau^X_D E^y\tau^X_D \leq G^X_D(x,y).
\end{equation*}
Therefore from Theorem \ref{general} we infer that
\begin{cor}
Let $p^X_D(t,\cdot,\cdot)$ be continuous for every $t>0,$ and let a finite measure
$\sigma=\nu^X-\nu^Y$ be nonnegative.
  Suppose that the Lebesgue
measure is absolutely continuous with respect to $\nu^X$. Then for
any bounded open set $D$ there exists a constant
$C=C(\sigma,D,\alpha,d)$ such that for $x\in D$,
\[G^Y_D(x,y)\leq C G^X_D(x,y),\quad (y)\text{ }a.s..\]
\end{cor}

Our next goal is  to reverse  the above estimate. We are not able
to do it under the above  assumptions but this will be done under
some additional assumptions through several steps. In the first
one  we take advantage of the following lemma which can be proved
similarly as Lemma 7 in \cite{R}.

\begin{lem}
Let $\sigma=\nu^X-\nu^Y$ be a nonnegative finite measure. Suppose
that $G^X_D(x,\cdot)$ and $G^Y_D(x,\cdot)$ are continuous then
$$G^X_D(x,y)\leq G^Y_D(x,y)+E^x[\tau^X_D>T;G^X_D(X_T,y)], $$
where $T$ is defined by (\ref{defT}).
\end{lem}
This lemma can be rewritten in the way which is more useful for further analysis.
\begin{cor}\label{Gformula}
Suppose that $\sigma=\nu^X-\nu^Y$ is a nonnegative finite measure,
$G^X_D(x,\cdot)$ and $G^Y_D(x,\cdot)$ are continuous. Then
$$G^X_D(x,y)\leq G^Y_D(x,y)+\int_D\int_{D-w} G^Y(x,w)G^X(w+z,y)\sigma(dz)dw.$$
\end{cor}
\begin{proof}
See the proof of Lemma 9 in \cite{R}.
\end{proof}

From now on  we assume that $X_t= \tilde{X}_t$ and  that the
measure $\sigma=\til{\nu}-\nu^Y$ is finite and absolutely
continuous. We will use the following notational convention:  in
the case when a measure $\mu$ is absolutely continuous we denote
its density by $\mu(x)$. That is $\sigma(x)$ is the density of
$\til{\nu}-\nu^Y$ Moreover we assume a particular  behavior of
$\sigma(x)$ near $0$, that is we suppose there exist $\varrho>0$
and $C$ such that
\begin{equation}\label{sigma}|\sigma(x)|\leq C|x|^{\varrho-d}, \quad |x|\leq 1. \end{equation}
In addition we assume that $\sigma(x)$ is bounded on $B^c(0,1)$,
which obviously is equivalent to boundness of $\nu^Y(x)$ on
$B^c(0,1)$.

For example the above conditions are satisfied by the L\'{e}vy measure
of the  relativistic process (see \cite{R}) and the L\'{e}vy measure of
the $\alpha$-stable process truncated to $B(0,1)$
($\nu^Y(x)=\textbf{1}_{B(0,1)}(x)\til{\nu}(x)$).

With these assumptions we have that the characteristic function of
$Y_t$ is integrable, so $p^Y(t,\cdot)$ is bounded and continuous.
Moreover, by  (\ref{ogrGestStab}) we get that for any
$\delta>0$,
 $$\til{p}(t,x)\leq C(\delta),\quad |x|\geq\delta.$$
Therefore from Lemma \ref{BoundDensity} we obtain that
the transition density of $Y_t$ also satisfies
$$p^Y(t,x) \leq C(\delta)\quad  |x|\geq \delta.$$
This property enables us  to  prove, similarly as for the Brownian
motion in \cite{CZ}, that $p^Y_D(t,x,\cdot)$ and $p^Y(t,\cdot,y)$
are continuous, moreover  $G^Y_D(x,\cdot)$ and $G_D^Y(\cdot,y)$
are continuous, too. Hence under the present assumptions, in all
claims of the results proved so far, we have that the  estimates
hold for every $y$ not for almost all.

Furthermore, we have that there exists a radius $r$ and a constant $c$ such that
$\til{\nu}(x) \le c \nu^Y(x)$ on $B(0,r)$. So, $\inf_{x\in
B(0,r)}\nu^Y(x)>0$. Therefore from Theorem \ref{IUThr:this}  we have
that for any bounded and connected Lipschitz domain the process $Y_t$
satisfies property A. That is we have the following corollary.
\begin{cor}\label{propertyAhold }
Let $\sigma(x)=\til{\nu}(x)-\nu^Y(x)$ be an integrable function
satisfying (\ref{sigma}).  Moreover let $\sigma$ be bounded on $B^c(0,1)$.
Then the property A holds for $Y_t$ and  any bounded connected Lipschitz domain.
Whereas if we assume that $\nu^Y\geq \til{\nu}$ then the property A
holds  for $Y_t$ and  any bounded open set.
\end{cor}

Let $D\subset \Rd$ be a  bounded Lipschitz domain with Lipschitz character
$(r_0,\lambda)$ (see \cite{J}, \cite{Bogdan1} for the
definitions). We need to introduce some additional notation
related to $D$. We assume that $D$ is a nonempty, open and bounded
set. We put $\ro=\frac{r_0}{\diam(D)}$ and
$\kappa=1/(2\sqrt{1+\lambda^2})$. The set $\{x\in
D:\delta_D(x)\geq r_0/2\}$ is nonempty. We choose one of its
elements and denote by $x_0=x_0(D)$. Besides we fix a point $x_1$
such that $|x_0-x_1|=r_0/4$. For any $x,y\in D$ let
$r=r(x,y)=\delta_D(x)\vee \delta_D(y)\vee |x-y|$. If $r\leq
r_0/32$ we put $A_{x,y}$ as a element of the following set
$$\mathcal{B}(x,y)=\{A\in D : B(A,\kappa r)\subset D\cap B(x,3r)\cap B(y,3r)\},$$
and if $r>r_0/32$ we set $A_{x,y}=x_1$.

For Lipschitz domains Jakubowski \cite{J} proved the following
theorem about estimates of the Green function for the isotropic
$\alpha$-stable process  in the case $d\ge 2$. If $d=1$, then analogous theorem is
true as well for $\alpha<1$ (see e.g. \cite{Bycz1}).

\begin{thr}\label{jakubowski:this}
Let $D$ be a bounded Lipschitz domain and $d>\alpha$. There is a
constant $C_1=C_1(d,\lambda,r_0,\diam(D),\alpha)$ such that for
every $x,y\in D$ we have
\[C^{-1}_1\frac{\til{\phi}_D(x)\til{\phi}_D(y)}{\til{\phi}_D^2(A_{x,y})}\norm{x-y}^{\alpha-d}\leq
\G{x}{y} \leq
C_1\frac{\til{\phi}_D(x)\til{\phi}_D(y)}{\til{\phi}_D^2(A_{x,y})}\norm{x-y}^{\alpha-d},\]
where $\til{\phi}_D(x)=\G{x}{x_0}\wedge
\mathscr{A}(d,\alpha)r^{\alpha-d}_0$.
\end{thr}

From the scaling property of the Green function for the isotropic
$\alpha$-stable process we have the following remark.
\begin{rem}
The constant $C_1$ depends on $r_0$ and $\diam(D)$ only by their
ratio $\ro$.
\end{rem}

Now, we recall estimates for the Green function of the isotropic
$\alpha$-stable process if $1=d\leq\alpha$. Their  proof  can
be found e.g. in \cite{Bycz1}.
\begin{thr}\label{GtilOne}
Let $d=1$ and $D$ be an open interval. Then we have on $D\times
D$,
\[\G{x}{y} \approx\left\{%
\begin{array}{ll}
    \log\left(\frac{(\delta_D(x)\delta_D(y))^{1/2}}{|x-y|}+1\right), & \hbox{$\alpha=1$,}\\
\min\left\{(\delta_D(x)\delta_D(y))^{(\alpha-1)/2},
\frac{(\delta_D(x)\delta_D(y))^{\alpha/2}}{|x-y|}\right\}, & \hbox{$1<\alpha$.} \\\end{array}%
\right.\]
\end{thr}

The consequence of Lemma 13 and 15 from \cite{J} is the following
lemma.
\begin{lem}\label{lemma3phi:this}
There are constants $\gamma=\gamma(d,\lambda,\alpha)<\alpha<d$ and
$C=C(d,\lambda,\alpha,\ro)$ such that for every $x,y,z,w\in D$ we
have
\[\frac{\til{\phi}_D^2(A_{x,y})}{\til{\phi}_D(A_{x,w})\til{\phi}_D(A_{z,y})}\leq C
\max\left\{1,\frac{|x-y|^{\gamma}}{|x-w|^{\gamma}},
\frac{|x-y|^{\gamma}}{|z-y|^{\gamma}},
\frac{|x-y|^{2\gamma}}{|x-w|^{\gamma}|z-y|^{\gamma}}\right\}.\]
\end{lem}
\begin{proof}
First, we assume that $|x-y|\leq|x-w|$. Then it can be proved using
similar methods as in  Lemma  13 of \cite{J} that
\begin{equation}\til{\phi}_D(A_{x,y})\leq C(d,\lambda,\alpha,\ro)\til{\phi}_D(A_{x,w}).\label{phi1:this} \end{equation}
Now, let $|x-w|\leq|x-y|$. Then from the proof of Lemma 15 in
\cite{J} we infer that
\begin{equation}\til{\phi}_D(A_{x,y})\leq C(d,\lambda,\alpha,\ro)
\frac{|x-y|^{\gamma}}{|x-w|^{\gamma}}
\til{\phi}_D(A_{x,w}),\label{phi2:this}\end{equation} for some
$0<\gamma<\alpha$. Combining (\ref{phi1:this}) and
(\ref{phi2:this}) ends the proof.
\end{proof}

\begin{lem}\label{calka:this}
Let $x\neq y\in D$, $-d<\varrho$ and $0<a,b$. Then there exists a
constant $C=C(d,a,b,\varrho)$ such that
\[\int_D\int_D
|y-z|^{a-d}|z-w|^{\varrho}|w-x|^{b-d}dzdw \leq C
\begin{cases} |x-y|^{a+\varrho+b},& a + \varrho+b<0,\\
1+\log\left(\frac{\diam(D)}{|x-y|}\right),&a + \varrho+b=0,\\
(\diam(D))^a\left(1+\log\left(\frac{\diam(D)}{|x-y|}\right)\right),&a=b=-\varrho,\\
(\diam(D))^{a+\varrho+b},& otherwise.\end{cases}\]
\end{lem}
\begin{proof}
By changing  variables: $u=\frac{z-y}{|x-y|}$ and
$v=\frac{w-x}{|x-y|}$ we get
\[\int_D\int_D
|y-z|^{a-d}|z-w|^{\varrho} |w-x|^{b-d}dzdw =
|x-y|^{a+b+\varrho}\int_{\frac{D-y}{|x-y|}}\int_{\frac{D-x}{|x-y|}}
|u|^{a-d}|v|^{b-d}|u-v-\mathbf{q}|^{\varrho}dudv,\]
where $\mathbf{q}=\frac{x-y}{|x-y|}$.\\
For $\varrho+a<0$ we have
\[\int_{\Rd}
|u|^{a-d}|u-v-\mathbf{q}|^{\varrho}du=C_{d,a,\varrho}|v+\mathbf{q}|^{a+\varrho},\]
and for $\varrho+a+b<0$,
\[\int_{\Rd}
|v|^{b-d}|v+\mathbf{q}|^{a+\varrho}dv=C_{d,a,b,\varrho},\] which
proves the first case. When $\varrho+a+b=0$, then we have
\begin{eqnarray*}\int_{\frac{D-x}{|x-y|}}
|v|^{b-d}|v+\mathbf{q}|^{a+\varrho}dv&\leq&
\int_{B(0,2)}|v|^{b-d}|v+\mathbf{q}|^{a+\varrho}dv+
2^{-\varrho-a}\int_{B(0,\diam(D)/|x-y|) \backslash
B(0,2)}|v|^{-d}dv\\
&=&C(d,a,b,\varrho)+C(d,a,\varrho)\left(\log\left(\frac{\diam(D)}{|x-y|}\right)-\log(2))\right)\vee
0\\
&\leq&C(d,a,b,\varrho)\left\{1+\log\left(\frac{\diam(D)}{|x-y|}\right)\right\}.
\end{eqnarray*}
If $0<\varrho+a+b<b$ then
\begin{eqnarray*}\int_{\frac{D-x}{|x-y|}}
|v|^{b-d}|v+\mathbf{q}|^{a+\varrho}dv&\leq&
\int_{B(0,2)}|v|^{b-d}|v+\mathbf{q}|^{a+\varrho}dv+
2^{-\varrho-a}\int_{B(0,\frac{\diam(D)}{|x-y|}) \backslash
B(0,2)}|v|^{\varrho+a+b-d}dv\\
&\leq&C(d,a,b,\varrho)\left\{1+\left(\frac{\diam(D)}{|x-y|}\right)^{\varrho+a+b}\right\}.
\end{eqnarray*}
The remaining cases can be proved in the same way.
\end{proof}

\begin{lem}\label{GPGkernel:this}Let $d>\alpha$. Suppose that there is a positive $\varrho$
and $c_1=c_1(\diam(D))$ such that $|\sigma(x)|\leq c_1
|x|^{\varrho-d}$ for $|x|\le \diam(D)$. Then there exists a
constant $C=C(d,\lambda,\ro,\alpha,\varrho)$ such that for all
$x,y\in D$,
\[ \int_D\int_D \G{y}{z} |\sigma(z-w)| \G{w}{x} dw dz \leq c_1 C\, (\diam(D))^{\zeta_1} |x-y|^{\zeta_2} \G{x}{y},\]
for some $\zeta_1\geq0$ and $\zeta_2>0$.
\end{lem}

\begin{proof}
From Theorem \ref{jakubowski:this} and Lemma 13 in \cite{J} we
obtain
\begin{eqnarray*}\frac{\G{x}{w} \G{z}{y}}{\G{x}{y}}&\approx&
\left(\frac{|x-y|}{|x-w||y-z|}\right)^{d-\alpha}
\frac{\til{\phi}_D(w)\til{\phi}_D(z)\til{\phi}_D^2(A_{x,y})}{\til{\phi}_D^2(A_{x,w})\til{\phi}_D^2(A_{z,y})}\\&\leq&
\left(\frac{|x-y|}{|x-w||y-z|}\right)^{d-\alpha}
\frac{\til{\phi}_D^2(A_{x,y})}{\til{\phi}_D(A_{x,w})\til{\phi}_D(A_{z,y})}.\end{eqnarray*}
Because $|\sigma(x)|\leq c_1 |x|^{\varrho-d}$ for $|x|\le
\diam(D)$ we get $|\sigma(w-z)|\le c_1 |w-z|^{\varrho-d}$ on
$D\times D$. So, from Lemma \ref{lemma3phi:this} it's enough to
prove that for some $\zeta_1\geq0$ and $\zeta_2>0$,
\[|x-y|^{d-\alpha+\rho_1+\rho_2}\int_D\int_D
|x-w|^{\alpha-\rho_1-d}|w-z|^{\varrho-d}|z-y|^{\alpha-\rho_2-d}dwdz
\leq C (\diam(D))^{\zeta_1}|x-y|^{\zeta_2},\] for some
$C=C(d,\rho_1,\rho_2,\varrho)$, where
$\rho_1,\rho_2\in\{0,\gamma\}$. Recall that $\gamma<\alpha$, hence
the above inequality  is a consequence of Lemma \ref{calka:this}.
\end{proof}

By inspecting the estimates from Theorem \ref{GtilOne} one can
check that the following remark is true.
\begin{rem}
In the case $d=1\leq\alpha$ the above lemma does not hold. This
is a reason why the proof below of Theorem \ref{main1} in the
one-dimensional  case for $\alpha \ge 1$ needs to employ some
other arguments then in the general case.
\end{rem}

\subsection{Proof of Theorem \ref{main1}}
Throughout this subsection we assume that $\sigma=\til{\nu}-\nu^Y$ is a finite
nonnegative absolutely continuous measure and its density satisfies
$$ \sigma(x)\leq C |x|^{\varrho-d},\quad |x|\le 1,$$
for
some positive $\varrho$. Then there is also a constant
$c=c(C,d,\alpha,\diam(D))$ such that $\sigma(x)\leq c
|x|^{\varrho-d}$ for $|x|\le \diam(D)$.
Let $D$ be a bounded connected Lipschitz domain. Then the property A holds
for $Y_t$ by Theorem \ref{IUThr:this}.

 The corollaries \ref{stableGreenUp} and \ref{Gformula} allow us to
write the following inequality
\begin{equation}\label{GRG}C_1^{-1}G^Y_D(x,y)\leq \G{x}{y}\leq G^Y_D(x,y)+C_1\til{R}_D(x,y),
\end{equation} where $\til{R}_D(x,y)=\int_D\int_D \G{x}{w}\sigma(w-z)\G{z}{y}dwdz$.

From Theorems \ref{jakubowski:this} and \ref{GtilOne} we obtain that for $|x-y|\geq\theta>0$
$$ \G{x}{y}\leq C(\theta)E^x\til{\tau}_D E^y\til{\tau}_D.$$
Hence, by the property A and Lemma \ref{momentyogr:this} we get
$$\G{x}{y}\leq C(\theta)G^Y_D(x,y),\quad |x-y|\geq\theta>0 $$
What remains it is to show that  $\til{R}_D(x,y)\leq \frac{1}{2
C_1}\G{x}{y}$ if $|x-y|$ is small enough. But for $d>\alpha$ this
is a consequence of Lemma \ref{GPGkernel:this}. This completes the
proof for $d>\alpha$.

Now, we  deal with  the case  $1=d\leq \alpha$. We  need to show
that $\G{x}{y}\leq C G^Y_D(x,y)$ if $|x-y|$ is small enough.
Recall that in this case $D$  is a bounded open interval.

\begin{lem}\label{Rfala:this}Let $d=1$.
Then there is a constant $C=C(\alpha,D,m)$ such that for any
$x,y\in D$,
\[\til{R}_D(x,y)\leq C\frac{(\delta_D(x)\delta_D(y))^{\alpha/2}}
{|x-y|^{1-\varrho\wedge 1}}.\]
\end{lem}
\begin{proof}
From Theorem \ref{GtilOne} it is easy to see that
\begin{equation}\label{GOneGora}
\G{x}{y}\leq C \frac{(\delta_D(x)\delta_D(y))^{\alpha/2}}{|x-y|}.
\end{equation}
Hence, for $\varrho<1$ we can prove in the same way as in  Lemma 8 in \cite{R}
that \begin{equation}\label{Rfalaproof1}\int_D
\widetilde{G}_D(x,w)\frac{dw}{|w-y|^{1-\rho}}\leq C
\frac{(\delta_D(x))^{\alpha/2}}{|x-y|^{1-\rho}}.\end{equation}
From the above
\begin{eqnarray*}
\int_D \widetilde{G}_D(x,w)\sigma(z-w)dw&\leq& C \int_D
\widetilde{G}_D(x,w)\frac{dw}{|w-z|^{1-\varrho}}\leq C
\frac{\delta_D(x)^{\alpha/2}}{|x-z|^{1-\varrho}}.
\end{eqnarray*}
If $\varrho\ge 1$ then $\sigma$ is bounded and one knows that
$E^x\til{\tau}_D\approx (\delta_D(x))^{\alpha/2}$ , so
$$\int_D \widetilde{G}_D(x,w)\sigma(z-w)dw\le C E^x\til{\tau}_D\le
c \delta_D(x)^{\alpha/2}.$$ Now, we use symmetry of the Green
function and the inequality \ref{Rfalaproof1} again to get
\[\widetilde{R}_D(x,y)\leq C \frac{(\delta_D(x)\delta_D(y))^{\alpha/2}}{|x-y|^{1-\varrho\wedge1}}.\]
\end{proof}

Finally, we are able to prove the lower bound of the Green function
for $1=d\le \alpha$.
\begin{prop}
\label{GDolOne}
Let $D$ be a bounded and open interval. Let $\alpha\ge 1$. Then
there exists a constant $C=C(m,d,\alpha,D)$ such that for any
$x,y\in D$,
\[\G{x}{y}\leq C G^Y_D(x,y).\]
\end{prop}

\begin{proof}
Note that we only need to consider the case $|x-y|\le \theta$  for
some sufficiently small $\theta>0$. First, we assume that $\D\leq
|x-y|^2$. By Theorem \ref{GtilOne} this implies that
\[\frac{(\D)^{\alpha/2}}{|x-y|}\leq C \G{x}{y}. \]
Then apply Lemma \ref{Rfala:this} to obtain
\[\widetilde{R}_{D}(x,y)\leq C |x-y|^{\varrho\wedge 1}\widetilde{G}_D(x,y),\]
for some constant $C$. So from (\ref{GRG}) it follows that
\begin{equation}\label{Gfunction1}
\G{x}{y} \leq G^Y_D(x,y)+\widetilde{C}|x-y|^{\varrho\wedge 1}\G{x}{y}.
\end{equation}
  By
the estimates of $\widetilde{p}_D(t,x,y)$ (Remark
\ref{gestzabityStable}) we have
\begin{equation}
\int^{\infty}_{t_0}\widetilde{p}_D(t,x,y)dt \leq C
t^{-1-1/\alpha}_{0} (\D)^{\alpha/2}. \label{integral10}
\end{equation}
 Next, from Lemma \ref{gestzabity1:this} for $X=\til{X}$ we have
\begin{equation}\label{ogrpDfala:this}
\widetilde{p}_D(t,x,y)\leq
p^Y_D(t,x,y)+\widetilde{p}(t,x,y)-e^{-2mt}p^Y(t,x,y),
\end{equation}
so integrating over $[0,t_0]$, where $t_0=(\D)^{\alpha/6}\le 1 $,
using Lemma \ref{potlematuzup:this}, and combining with
(\ref{integral10}) we obtain
\begin{eqnarray}
\G{x}{y}&=& \int_{0}^{t_0}\widetilde{p}_D(t,x,y)dt+\int^{\infty}_{t_0}\widetilde{p}_D(t,x,y)dt\notag\\
&\leq& G^Y_D(x,y)+
\int^{t_0}_0(\widetilde{p}(t,x,y)-e^{-2mt}p^Y(t,x,y))dt+C
t^{-1-1/\alpha}_{0} (\D)^{\alpha/2}\notag\\
&\leq& G^Y_D(x,y) + c t^{2-1/\alpha}_0 +C
t^{-1-1/\alpha}_{0} (\D)^{\alpha/2}\notag\\
&=& G^Y_D(x,y) + c (\D)^{\frac{2\alpha-1}6}.\label{Gfala3:this}
\end{eqnarray}
Now assume that $|x-y|^2\leq\D$ and  take into account that in
this case $\G{x}{y}\ge C (\D)^{(\alpha-1)/2}$, so  we can rewrite
(\ref{Gfala3:this}) as
\begin{equation}\label{deltaEstimate}
\G{x}{y}\leq G^Y_D(x,y)+ c (\D)^{\rho}\G{x}{y},
\end{equation}
where ${\rho}=\frac {2-\alpha}6>0$.
Observe that (\ref{deltaEstimate}) in the  case $|x-y|^2\leq\D\le
\theta$, and  (\ref{Gfunction1}) in the  case $\D\le |x-y|^2\le
\theta$ for $\theta$ sufficiently small provide the conclusion.
 From the remaining cases $\D\ge  \theta$ or $|x-y|^2\ge \theta$
 only the first needs to be considered and  can be
handled in a very simple way. Indeed,  in this situation
$$(\D)^{\frac{2\alpha-1}6}\le (\D)^{\frac{\alpha}2}\theta^{-\frac{\alpha+1}6}\le
C \theta^{-\frac{\alpha+1}6}G^Y_D(x,y),$$
where the last step follows from the fact that $Y_t$ has the
property A and Lemma \ref{momentyogr:this}. Hence the conclusion holds by (\ref{Gfala3:this}).
This completes the proof.

\end{proof}

\subsection{Case $\nu^Y\geq \til{\nu}$ }

Throughout this subsection we assume that $\nu^Y\geq \til{\nu}$
and in addition let $D$ be a bounded Lipschitz domain.  Note that
in this case by the result  of Sztonyk \cite{Sztonyk} the process
$Y$ does not hit the boundary on exiting $D$, so if $u$ is regular
harmonic on $D$ with respect to the process $Y$ then
\begin{equation}\label{boundry}
u(x)=E^xu(Y_{\tau_{D}})=\int_{D^c}u(z)P^Y_{D}(x,z)dz, \quad x\in D.
\end{equation}
 The aim of
this section is to prove that the Green functions are comparable,
first for $D$ with small diameter and then for arbitrary bounded
Lipschitz domains.
The result for $D$ of small diameter allows us
to prove a version of the Boundary Harnack Principle under the following
assumptions :
\begin{description}
    \item[G1] $\nu^Y(x)\geq \til{\nu}(x)$ for $x\in\Rd \setminus\{0\}$,
    \item[G2] for some $R>0$ there are constants $c_1(R)$ and $\gamma$ such that
    $$ |\sigma(x)|=|\til{\nu}(x)-\nu^Y(x)|\leq c_1
    |x|^{\varrho-d}\quad\text{for } |x|\leq R,$$
    \item[G3] there is a constant $c_2=c_2(R)$ such that
$$ \nu^Y(x)\leq c_2 \nu^Y(y)\quad \text{for any }x,y\in \Rd\text{ such
that }|x-y|\leq R/2 \text{ and } |x|,|y|\geq R/2.$$
\end{description}
Then after establishing BHP we show that we can remove the assumption about the diameter
of the set $D$.

We start with the iteration of the inequality from Corollary
(\ref{Gformula}) to obtain for $G^Y_D(x,\cdot)$ continuous,
\begin{eqnarray}\label{GPGIter}
G^Y_D(x,y)&\leq& \til{G}_D(x,y)+
\sum^{n}_{k=1} [(H_D^{\sigma})^k \til{G}_D(\cdot,y)](x) +
[(H_D^{\sigma})^{n+1} G^Y_D(\cdot,y)](x),
\end{eqnarray}
where  $H^\sigma_D: L^1(D)\rightarrow L^1(D)$ is given by
$$ [H_D^{\sigma} f(\cdot)](x)=\int_D\int_D\G{x}{w}|\sigma(w-z)|f(z)dw dz.$$


We now prove comparability of Green functions for sets of small
diameter.  Note that the constant $C$ in the conclusion of the
following Proposition depends on $D$ through $\ro$ and $\lambda$.
This feature is crucial for our future applications.

\begin{prop}\label{Greensmall:this}Let $d>\alpha$.
Let $D$ be a Lipschitz domain and $G^Y_D(x,\cdot)$ be continuous
and $\nu^Y$ satisfies \textbf{G1} and \textbf{G2}, then there
exist constants $R_0=R_0(d,\alpha,\lambda,\ro,\sigma)\leq R$ and
$C=C(R_0)$  which has the following property: if $\diam(D)\leq
R_0$ then
\[C^{-1}\til{G}_{D}(x,y)\leq G^Y_{D}(x,y)\leq C \til{G}_{D}(x,y),\quad x,y\in D.\]
\end{prop}

\begin{proof}
 If $\diam(D)\leq R$ by
Lemma \ref{GPGkernel:this} we get that
\begin{equation*}\label{GPlineGkernel:this}
[H^{\sigma}_D\G{\cdot}{y}](x)\leq C_1
\diam(D)^{\zeta} \G{x}{y},
\end{equation*}
for some constant $C_1=C_1(d,\alpha,\lambda,\ro,\sigma)$ and
$\zeta>0$. Iterating the above inequality we obtain that
$[(H^{\sigma}_D)^k\G{\cdot}{y}](x)$ is bounded by
\[(C_1\diam(D)^{\zeta})^k\G{x}{y}.\]
 Setting
\[R_0=\frac{1}{2}C_1^{-1/\zeta}\wedge R \] we obtain for  $\diam(D)\leq R_0$ that
\begin{equation}\label{Greensmallproof1}[H^{\sigma}_D\G{\cdot}{y}](x)\leq \theta \G{x}{y},\end{equation}
for some $\theta\leq 1/2$.

Next, we show that for any $x\neq y \in D$
$$\lim_{n\rightarrow\infty}[(H_D^{\sigma})^{n} G^Y_D(\cdot,y)](x)=0. $$
Indeed, let us  observe that for a positive $f\in L^1(D)$ we have from
(\ref{Greensmallproof1}) that
\begin{eqnarray*}
[(H_D^{\sigma})^2f](x)&=& \int_D\int_D\int_D\int_D
\G{x}{u}|\sigma(u-v)| \G{v}{w}|\sigma(w-z)|f(z) dz dw dv du\\
&=&\int_D\int_D
[(H_D^{\sigma})\til{G}_D(\cdot,w)](x)|\sigma(w-z)|f(z) dz dw\\
&\leq & \theta \int_D\int_D \til{G}_D(x,w)|\sigma(w-z)|f(z) dz
dw\\
&=&\theta[H_D^{\sigma}f](x).
\end{eqnarray*}
Iterating  we obtain $[(H_D^{\sigma})^{n+1} G^Y_D(\cdot,y)](x)\leq
\theta^n [(H_D^{\sigma})G^Y_D(\cdot,y)](x)$.
 So it is enough
to prove that $[(H_D^{\sigma})G^Y_D(\cdot,y)](x)$ is finite. But
from Lemma \ref{potentialcomp} we obtain that there is a constant
$C$ such that $G^Y(x,y)\leq C\til{U}(x-y)$. Hence by Lemma
\ref{calka:this} we get
$$[(H_D^{\sigma})G^Y_D(\cdot,y)](x)\leq C \int_D\int_D \til{U}(x-w)|\sigma(w-z)|\til{U}(z-y)dwdz<\infty.$$

Finally, we infer from (\ref{GPGIter}) that if $\diam(D)\leq R_0$
then
\[G^Y_{D}(x,y)\leq \frac{\theta}{1-\theta}\G{x}{y},\]
which together with Corollary \ref{stableGreenUp} ends the proof.
\end{proof}

\begin{rem}
The constant $C(R_0)$ in the above theorem converges to $1$ if $\diam(D)$ converges to $0$.
\end{rem}

The next result shows that the Poisson kernels for $D$ are
comparable  under the assumptions of the preceding result. This in
consequence provides necessary tools to establish BHP, which
is employed to show comparability of Green functions for sets  of arbitrary
finite diameter.

\begin{prop}\label{Poissonsmall:this}Let $d>\alpha$ and $D$ be a
bounded Lipschitz domain. Assume that $\nu^Y$ satisfies
assumptions \textbf{G1} and \textbf{G2} and is bounded on
$B^c(0,R)$. There exist constants
$R_0=R_0(d,\alpha,\lambda,\ro,\sigma)\leq R/2$ and $C=C(R_0)$
which satisfy for $D$ such that $\diam(D) \leq R_0$
\[C^{-1}\til{P}_{D}(x,z)\leq P^Y_{D}(x,z)\leq C \til{P}_{D}(x,z),\]
for any $x\in D$ and $z\in \overline{D}^c: \delta_{D}(z)\leq R_0$.
Moreover, if we suppose that $\nu^Y$ satisfies assumption
\textbf{G3} with $R=2\,R_0$, then there exists a constant $C(R_0)$
such that
\[C^{-1}\nu^Y(z-x)E^x\til{\tau}_{D}\leq P^Y_{D}(x,z)\leq C\nu^Y(z-x)E^x\til{\tau}_{D},\]
for $x\in D$ and $z\in \overline{D}^c:
\delta_{D}(z)>R_0$.
\end{prop}

\begin{proof}
By Proposition \ref{Greensmall:this} there are  constants
$\overline{R}_0\leq R/2$ and $C_1(\overline{R}_0)$ such that
\[C^{-1}_1\til{G}_{D}(x,y)\leq G^Y_{D}(x,y)\leq C_1\til{G}_{D}(x,y),\]
for $D$ with $\diam(D)\leq\overline{R}_0$. Next, from Theorem 1 in \cite{IW} we have the
following formula
\[P^Y_{D}(x,z)=\int_{D}\nu^Y(z-y)G^Y_{D}(x,y)dy.\]
But $|\sigma(w)|\leq c_1|w|^{-d+\varrho}= c_1
\mathscr{A}(-\alpha,d)^{-1}\til{\nu}(w)|w|^{\varrho+\alpha}$. So
for $z\in \overline{D}^c: \delta_{D}(z)\leq \overline{R}_0$ we
have
\[|\sigma(z-y)|\leq c_1 \mathscr{A}(-\alpha,d)^{-1}
(2\overline{R}_0)^{\varrho+\alpha}\til{\nu}(z-y).\] Hence, we put
$R_0=\overline{R}_0\wedge 1/2\left(\frac{\mathscr{A}(-\alpha,d)}
{2c_1}\right)^{1/(\alpha+\varrho)}$ and then
$$|\sigma(z-y)|\leq \frac{1}{2}\til{\nu}(x). $$
By the above inequality we obtain
\begin{eqnarray*}
P^Y_{D}(x,z)&\leq& C_1 \int_{D}\nu^Y(z-y)\til{G}_{D}(x,y)dy  \\
&= & C_1 \left(\int_{D}\til{\nu}(z-y)\til{G}_{D}(x,y)dy+
\int_{D}\sigma(z-y)\til{G}_{D}(x,y)dy \right)\\
&\leq& C_1 \til{P}_{D}(x,y) + C_1\int_{D}|\sigma(z-y)|\til{G}_{D}(x,y)dy\\
&\leq& \frac{3}{2}C_1 \til{P}_{D}(x,y),
\end{eqnarray*}
and
\begin{eqnarray*}
P^Y_{D}(x,z)&\geq& C^{-1}_1 \int_{D}\nu^Y(z-y)\til{G}_{D}(x,y)dy  \\
&\geq& C^{-1}_1 \til{P}_{D}(x,y) - C^{-1}_1\int_{D}|\sigma(z-y)|\til{G}_{D}(x,y)dy\\
&\geq& \frac{C^{-1}_1}{2} \til{P}_{D}(x,y),
\end{eqnarray*}
which ends the proof of the first claim of the theorem.

Now, suppose that there is a constant $c=c(R_0)$ such that
$\nu^Y(x)\leq c \nu^Y(y)$ for all $|x|,|y|\geq R_0$ such that
$|x-y|\le R_0$. Assume that  $z\in \overline{D}^c:
\delta_{D}(z)>R_0$. For $x,y\in D$ we have
\[|x-z|\geq \delta_{D}(z)\geq R_0\text{ and of course
} |x-y|\leq\diam(D)\leq R_0.\] Hence, we get
\begin{eqnarray*}
P^Y_{D}(x,z)&\leq&cC_1\nu^Y(x-z)\int_{D}\til{G}_{D}(x,y)dy\\
 &=&cC_1\nu^Y(x-z)E^x\til{\tau}_{D}.
\end{eqnarray*}
Similarly the lower bound is
\[P^Y_{D}(x,z)\geq (cC_1)^{-1}\nu^Y(x-z)E^x\til{\tau}_{D}.\]
\end{proof}


\begin{thr}(Boundary Harnack Principle-BHP)\label{BHPsmall:this}
Let $d>\alpha$ and $D$ be a bounded Lipschitz domain. Suppose that
$\nu^Y$ satisfies \textbf{G1}-\textbf{G3}. Let $Z\in
\partial D$. Then there exists a constant $\rho_0=\rho_0(D)$ such
that for any $\rho\in(0,\rho_0]$ and  two functions $u$ and $v$
which are nonnegative in $\Rd$ and positive, regular harmonic in
$D\cap B(Z,\rho)$. If $u$ and $v$ vanish on $D^c\cap B(Z,\rho)$,
then for $x,y\in D\cap B(Z,\rho\beta)$
\[\frac{u(x)}{v(x)}\leq C \frac{u(y)}{v(y)},\]
for some constant $C=C(D,\alpha,\sigma)$ and $\beta(d,\lambda)\in (0,1)$.
\end{thr}

\begin{proof}
There is a constant $R_1=R_1(d,\lambda)\geq 1$ (see e.g.
\cite{Bogdan1})  such that for all $Z\in\partial D$ and
$r\in(0,r_0)$, there exists a Lipschitz domain $\Omega(r)$  with
the Lipschitz constant $\lambda\, R_1$ and the localization radius
$\diam(D) \ro /R_1$, having the property
$$D\cap B(Z,r/R_1)\subset \Omega(r)\subset D\cap B(Z,r).$$

The proof consists of  showing that there are constants
$C=C(D,\alpha,\sigma)$ and $\rho_0$ such that for  $\rho<\rho_0$ and $z\in \Omega(\rho)^c\cap
B^c(Z,\rho/R_1)$,
\begin{equation}\label{BHPproof1}
P^{Y}_{\Omega(\rho)}(x,z)\leq C
\frac{E^x\til{\tau}_{\Omega(\rho)}}{E^y\til{\tau}_{\Omega(\rho)}}
P^{Y}_{\Omega(\rho)}(y,z),
\end{equation}
where $x,y\in D\cap B(Z,\rho/(R_1 2))$. It is worth mentioning
that the constant $C$  is universal for all sets $\Omega(\rho),\
\rho\le \rho_0$. This would  give the conclusion with
$\beta=1/(2\, R_1)$ since by (\ref{boundry}) we have

\begin{eqnarray*}
u(x)=E^xu(Y_{\tau_{\Omega(\rho)}})&=&\int_{\Omega(\rho)^c}u(z)P^Y_{\Omega(\rho)}(x,z)dz\\
&=&\int_{\Omega(\rho)^c \backslash B(Z,\rho/R_1)}u(z)P^Y_{\Omega(\rho)}(x,z)dz\\
&\leq&C  \frac{E^x\til{\tau}_{\Omega(\rho)}}{E^y\til{\tau}_{\Omega(\rho)}}
\int_{\Omega(\rho)^c\backslash B(Z,\rho/R_1)}u(z)P^Y_{\Omega(\rho)}(y,z)dz\\
&= & C \frac{E^x\til{\tau}_{\Omega(\rho)}}{E^y\til{\tau}_{\Omega(\rho)}}
u(y),
\end{eqnarray*}
which would imply
\[\frac{u(x)}{u(y)}\frac{v(y)}{v(x)}\leq  C\frac{E^x\til{\tau}_{\Omega(\rho)}}{E^y\til{\tau}_{\Omega(\rho)}}
C\frac{E^y\til{\tau}_{\Omega(\rho)}}{E^x\til{\tau}_{\Omega(\rho)}}=C^2.\]

Now we prove (\ref{BHPproof1}). From Proposition
\ref{Poissonsmall:this} we obtain that there exists constant
$\rho_0<r_0(D)$ and $C_1=C_1(\rho_0)$ such that for any $\rho\leq
\rho_0$
\[C_1^{-1}\til{P}_{\Omega(\rho)}(x,z)\leq
P^Y_{\Omega(\rho)}(x,z) \leq C_1 \til{P}_{\Omega(\rho)}(x,z),
\]
if $\delta_{\Omega(\rho)}(z)\leq \rho_0$. Note that $C_1$ is universal for all $\Omega(\rho)$.

By Theorem 2 in \cite{J} we have that there is some $C_2=C_2(\alpha,d,\lambda,\ro)$
such that for any $x,y\in D$ and $z\in \overline{D}^c$
\[\til{P}_{\Omega(\rho)}(x,z)\leq C_2 \frac{E^x\til{\tau}_{\Omega(\rho)}}{E^y\til{\tau}_{\Omega(\rho)}}
\frac{\til{\phi}^2_{\Omega(\rho)}(A_{y,z'})}{\til{\phi}^2_{\Omega(\rho)}(A_{x,z'})}
\frac{|y-z|^{d-\alpha}}{|x-z|^{d-\alpha}}
\til{P}_{\Omega(\rho)}(y,z),\] where $z'\in \{A\in D:
B(A,\kappa\delta_{\Omega(\rho)}(z))\subset D\cap
B(S,\delta_{\Omega(\rho)}(z)) \}$ if $\delta_{\Omega(\rho)}(z)\leq
r_0/32$ and $z'=x_1$ if $\delta_{\Omega(\rho)}(z)>r_0/32$ for $S$
such that $|z-S|=\delta_{\Omega(\rho)}(z)$. If $x,y\in D\cap
B(Z,\rho/(R_1 2))$ and $z\in \Omega(\rho)^c\cap B^c(Z,\rho/R_1)$
then
\[\frac{|y-z|}{|x-z|}\leq\frac{|x-z|+|x-y|}{|x-z|}\leq (1+\frac{\rho/R_1}{\rho/(2R_1)})=3.\]
Now, suppose that $\delta_{\Omega(\rho)}(z)\leq \rho/32$ then we obtain
\[|x-z'|\geq |x-z|-|z-z'| \geq |x-z|-|z-S|-|z'-S|\geq
\frac{\rho}{2}-2\delta_{\Omega(\rho)}(z)\geq
\frac{7}{16}\rho>\frac{r_0}{32},\] while if
$\delta_{\Omega(\rho)}(z)> \rho/32$ then $z'=x_1$, so
$\delta_{\Omega(\rho)}(z')\geq r_0/4$. Therefore
$A_{x,z'}=x_1=A_{y,z'}$ and of course
$\frac{\til{\phi}_{\Omega(\rho)}(A_{y,z'})}{\til{\phi}_{\Omega(\rho)}(A_{x,z'})}=1$.
Hence for $x,y\in D\cap B(Z,\rho/(R_1 2))$  and $z\in
\Omega(\rho)^c\cap B^c(Z,\rho/R_1)$ such that
$\delta_{\Omega(\rho)}(z)\leq \rho_0$ we get
$$ P^{Y}_{\Omega(\rho)}(x,z)\leq C^{2}_1 C_2 3^{d-\alpha}
\frac{E^x\til{\tau}_{\Omega(\rho)}}{E^y\til{\tau}_{\Omega(\rho)}} P^{Y}_{\Omega(\rho)}(y,z).$$

Next, observe that \textbf{G1}-\textbf{G3} imply that for $r\le R$
there is a constant  $c=c(r)$ such that $\nu^Y(x)\leq c \nu^Y(y)$
for all $x$ and $y$ such that $|x-y|\leq r$ and $|x|,|y|\geq r$.
Hence for $\delta_{\Omega(\rho)}(z)\geq \rho_0$ we have
\[P^Y_{\Omega(\rho)}(x,z) \leq C_3(\rho_0) \nu^Y(z-x)E^x\til{\tau}_{\Omega(\rho)}\leq C_3(\rho_0) c(\rho_0)
\nu^Y(z-y)E^x\til{\tau}_{\Omega(\rho)}\leq
cC^2_3\frac{E^x\til{\tau}_{\Omega(\rho)}}{E^y\til{\tau}_{\Omega(\rho)}}P^Y_{\Omega(\rho)}(y,z).\]

This completes the proof of (\ref{BHPproof1}) and hence the theorem.
\end{proof}
For regular harmonic functions, which vanish on $D^c$ we infer the
following remark.
\begin{rem}\label{BHPGreen}
Suppose $\nu^Y$ satisfies \textbf{G1}, \textbf{G2} and is bounded
on $B^c(0,R)$. Let $Z \in \partial  D$. Then there exists a
constant $\rho_0=\rho_0(D)$ such that for any $\rho\in(0,\rho_0]$
and  two functions $u$ and $v$ which are nonnegative in $\Rd$ and
positive, regular harmonic in $D\cap B(Z,\rho)$. If $u$ and $v$
vanish on $D^c$, then for $x,y\in D\cap B(Z,\rho\beta)$
\[\frac{u(x)}{v(x)}\leq C \frac{u(y)}{v(y)},\]
for some constant $C=C(D,\alpha,\sigma)$ and $\beta(d,\lambda)\in
(0,1)$.
\end{rem}

\begin{thr} \label{GreenBig}Let $d>\alpha$ and $D$ be a bounded Lipschitz domain.
Assume that $\nu^Y$ satisfies assumptions \textbf{G1}, \textbf{G2}
and is bounded on $B^c(0,R)$. Then for $x,y\in D$ we have
\[C^{-1}\G{x}{y}\leq G^Y_D(x,y)\leq C \G{x}{y},\]
for some constant $C=C(d,\lambda,\ro,\sigma)$.
\end{thr}

\begin{proof}
   Observe that for $|x-y|\leq N
   (\delta_D(x)\wedge\delta_D(y))$,
  \begin{equation*}
G^Y_D(x,y)\geq G^Y_{B(x,\delta_D(x)\wedge\delta_D(y)\wedge R_0(D))}(x,y)
  \geq C \til{G}_{B(x,\delta_D(x)\wedge\delta_D(y)\wedge
  R_0)}(x,y),
\end{equation*}
where $R_0$ is such that
  $G^Y_{B(0,R_0)}(x,y)\approx \til{G}_{B(0,R_0)}(x,y)$ (such $R_0$
  exists from Proposition \ref{Greensmall:this}). Next, it is easy to see from
Theorem 3.4 in \cite{Kulczycki} that
\begin{equation}\label{GreenBigproof1}
c(N)|x-y|^{\alpha-d}\leq
\til{G}_{B(x,\delta_D(x)\wedge\delta_D(y)\wedge
  R_0)}(x,y)\leq C G^Y_D(x,y).
\end{equation}
From Lemma \ref{potentialcomp} we have
\begin{equation}\label{GreenBigproof2}
G^Y_D(x,y)\leq U^Y(x-y)\leq C \til{U}(x-y) = C |x-y|^{\alpha-d}.
\end{equation}

We define similarly as in Theorem \ref{jakubowski:this} the
truncated Green function for $Y_t$ by
$$\phi^Y_D(x)=G^Y_D(x_1,y)\wedge \mathscr{A}(d,\alpha)r^{d+\alpha}_0.$$
Using Remark  \ref{BHPGreen}  we can repeat the arguments from
Lemma 17 in \cite{J} to show that
  \[\phi^Y_D(x)\approx E^x\tau^Y_D.\]
Next, by Lemma \ref{momentyogr:this} we get
    \[E^x \tau^Y_D\approx E^x \til{\tau}_D. \]
 Therefore
\begin{equation}\label{GreenBigproof3}
  \phi^Y_D(x)\approx \til{\phi}_D(x).
\end{equation}
  By the above and (\ref{GreenBigproof2}) we infer that there is a constant $r$ such that
$\phi^Y_D(x)=G^Y_D(x,x_0)$ for $x\in D\cap B^c(x_0,r)$. Hence by
Harnack's inequality for $\alpha$-stable harmonic functions we
obtain that for $x,y\in D\cap B^c(x_0,r)$  such that $|x-y|\leq N
(\delta_D(x)\wedge\delta_D(y)),$
\begin{equation}\label{GreenBigproof4}
G^Y_D(x,x_0)=\phi^Y_D(x)\approx\til{\phi}_D(x)\leq
  C(N)\, \til{\phi}_D(y)\approx\phi^Y_D(y)=G^Y_D(y,x_0).
\end{equation}

Using  BHP for $Y_t$ (Remark \ref{BHPGreen}), and taking into
account (\ref{GreenBigproof1}), (\ref{GreenBigproof2}) and
(\ref{GreenBigproof4}) we can prove a version of Theorem
\ref{jakubowski:this} with $G^Y_D$ instead of $\til{G}_D$ (see the
proof of Theorem 1 in \cite{J}), that is
$$C^{-1}_1\frac{\phi^Y_D(x)\phi^Y_D(y)}{(\phi^Y_D(A_{x,y}))^2}\norm{x-y}^{\alpha-d}\leq
G^Y_D(x,y) \leq
C_1\frac{\phi^Y_D(x)\phi^Y_D(y)}{(\phi^Y_D(A_{x,y}))^2}\norm{x-y}^{\alpha-d}.
$$ Applying  (\ref{GreenBigproof3}) and then comparing the above
estimate  with the bound from Theorem \ref{jakubowski:this} we get
the conclusion.
\end{proof}

\subsection{Proof of Theorem 1.2}
Let $d>\alpha$ and $D$ be a connected Lipschitz domain. Suppose
that $|\sigma(x)|\leq c_3 |x|^{-d+\varrho}$ for $|x|\leq 1$, where
$\varrho>0$ and $\nu^Y(x)$ is bounded on $B^c(0,1)$. Then the
property A holds for $Y_t$ from Corollary \ref{propertyAhold }.

Let $\{Z_t\}$ be a L\'{e}vy process with the L\'{e}vy measure,
which density is equal to $\nu(x)\vee \til{\nu}(x)$. Then of
course the process $Z_t$ and the set $D$ satisfies the assumptions
of Theorem \ref{GreenBig}. So, we obtain that there is a constant
$C_1$ such that
\begin{equation}\label{TH1.2proof1}C_1^{-1} \G{x}{y}\leq G^Z(x,y)\leq C_1 \G{x}{y}.\end{equation}
Therefore we have that
\begin{equation}\label{TH1.2proof2}C^{-1}_2\frac{\til{\phi}_D(x)\til{\phi}_D(y)}{(\til{\phi}_D(A_{x,y}))^2}
\norm{x-y}^{\alpha-d}\leq G^Z_D(x,y) \leq
C_2\frac{\til{\phi}_D(x)\til{\phi}_D(y)}{(\til{\phi}_D(A_{x,y}))^2}\norm{x-y}^{\alpha-d}.\end{equation}
Moreover, the property A holds for $Y_t$, that is
\begin{equation}\label{TH1.2proof3}C \, E^x\tau^Y_D E^y\tau^Y_D \leq G^Y_D(x,y).\end{equation}
Having (\ref{TH1.2proof2}) and (\ref{TH1.2proof3}) hold,  we can repeat the
proof of  Theorem \ref{main1} for $d>\alpha$. Hence there
exists a constant $C_3$ which satisfies
\begin{equation}\label{TH1.2proof4}
C^{-1}_3G^Y(x,y)\leq G^Z(x,y)\leq C_3 G^Y(x,y).
\end{equation}
Combining (\ref{TH1.2proof1}) and  (\ref{TH1.2proof4}) give us
\[C^{-1}\til{G}_{D}(x,y)\leq G^Y_{D}(x,y)\leq C \til{G}_{D}(x,y),\]
which completes the proof.

\end{document}